

\hsize=14cm
\vsize=20cm
\hoffset=1cm
\voffset=2cm
\footline={\hss{\vbox to 2.5cm{\vfil\hbox{\rm\folio}}}\hss}

\input amssym.def
\input amssym.tex

\font\csc=cmcsc10
\font\title=cmr12 at 14pt

\font\teneusm=eusm10    
\font\seveneusm=eusm7  
\font\fiveeusm=eusm5    
\newfam\eusmfam
\def\eusm{\fam\eusmfam\teneusm}
\textfont\eusmfam=\teneusm 
\scriptfont\eusmfam=\seveneusm
\scriptscriptfont\eusmfam=\fiveeusm

\font\teneufb=eufb10    
\font\seveneufb=eufb7  
\font\fiveeufb=eufb5    
\newfam\eufbfam

\textfont\eufbfam=\teneufb 
\scriptfont\eufbfam=\seveneufb
\scriptscriptfont\eufbfam=\fiveeufb


\def\Re{{\rm Re}\,}

\def\e#1{{\eusm #1}}

\centerline{\title A Multiple Sum Involving the M\"obius Function}
\vskip 1cm
\centerline{\csc Yoichi Motohashi}
\vskip 1cm
\noindent
{\csc Abstract.}\quad 
We consider a multiple arithmetical sum involving the M\"obius function
which despite its elementary appearance is in fact of a highly intriguing
nature. We establish an asymptotic formula for the quadruple case that raises
the first genuinely non-trivial situation. This is a rework of an
old unpublished note of ours.
\bigskip
\noindent
{\bf 1.\ Introduction.}\quad The aim of the present article is to discuss the
asymptotics of the quantity
$$
\e{M}_k(z)=\sum_{d_1\le z}\cdots
\sum_{d_{2k}\le z} 
{\mu(d_1)\cdots\mu(d_{2k})\over[d_1,\cdots,d_{2k}]},\quad k\ge1,\eqno(1.1)
$$
as $z$ tends to infinity, where $\mu$ is the M\"obius function and 
$[d_1,\cdots,d_{2k}]$ the least common multiple of positive
integers $d_1,\ldots, d_{2k}$. We have the relation
$$
\e{M}_k(z)=\lim_{N\to\infty}{1\over N}\sum_{n\le N}\left(\sum_{d|n, d\le
z}\mu(d)\right)^{2k}.\eqno(1.2)
$$
Thus our problem is pertinent to the extremal behaviour of the truncated
sum of the M\"obius function over divisors, and somewhat remotely to the
Selberg sieve (see the concluding remark). The case $k=1$ is treated in
[2]. The case $k=2$ is already quite involved and discussed in [10].
The result and the outline of the argument there have been shown on a few
occasions, first at the Problem Session of the Amalfi International Symposium on
Analytic Number Theory, September 1989. There is, however, a minor error in
[10], as is to be indicated below.
\medskip
We have reworked [10] because of its apparent relation with the recent important
challenge [5] by Goldston and Y{\i}ld{\i}r{\i}m to the problem of finding small
gaps between consecutive primes. Their argument depends on their previous work
[4] which deals with higher correlations of short sums of
$$
\Lambda_z(n)=\sum_{d|n, d\le z}\mu(d)\log(z/d),\eqno(1.3)
$$
an approximation to the von Mangold function. In their discussion it is
needed, among other things, to study the asymptotic behaviour of the sum
$$
\sum_{n\le N}\Lambda_z(n+j_1)\Lambda_z(n+j_2)\cdots\Lambda_z(n+j_r),\eqno(1.4)
$$
where $j_1,\ldots,j_r$ are arbitrary non-negative integers. Expanding this via
the definition $(1.3)$, we are led to an expression closely resembles to $(1.1)$.
Goldston and Y{\i}ld{\i}r{\i}m applied to this expression an argument
essentially the same as that of [10], apparently without being aware of our old
unpublished work. 
\par
It appears to us, however, that their problem is less delicate than ours,
as far as the handling of the relevant residue calculus is concerned. The
factor $\log(z/d)$ makes their expression smoother than ours. Being translated
into our situation, this is equivalent to having $(s_1\cdots s_{2k})^2$ in
place of the denominator $s_1\cdots s_{2k}$ in $(1.5)$ below. Hence, both the
convergence and the estimation issues are less troublesome with $(1.4)$,
although the arithmetical issue can be highly involved when $j_1,\ldots,j_r$
are arbitrary.
\medskip
The argument of [10] starts with the following integral expression: For
non-integral $z$
$$
\e{M}_k(z)={1\phantom{{}^{2k}}\over(2\pi i)^{2k}}\int\cdots\int
M(s_1,\cdots,s_{2k})z^{s_1+\cdots+s_{2k}}{ds_1\cdots
ds_{2k}\over s_1\cdots s_{2k}}\eqno(1.5)
$$
with
$$
M(s_1,\cdots,s_{2k})=\sum_{d_1=1}^\infty\cdots
\sum_{d_{2k}=1}^\infty {\mu(d_1)\cdots\mu(d_{2k})\over
[d_1,\cdots,d_{2k}]d_1^{s_1}\cdots d_{2k}^{s_{2k}}},\eqno(1.6)
$$
where all integral are over vertical lines placed in the right half plane.
Of course this is not a fully correct expression. We need to use, instead,
a truncated version of Perron's formula, and the vertical
segments over which the integrations are performed should be placed in a
well-poised way, as we shall show later.
\par
We have, for
$\Re s_j>0$ ($j=1,\ldots,2k$), the Euler product expansion
$$
M(s_1,\cdots,s_{2k})=\prod_p\left(1-{1\over p}+{1\over p}\prod_{j=1}^{2k}
\left(1-{1\over p^{s_j}}\right)\right).\eqno(1.7)
$$
Thus
$$
M(s_1,\cdots,s_{2k})={\prod\zeta(1+s_{\lambda_1}+\cdots+s_{\lambda_{2a}})
\over\prod\zeta(1+s_{\tau_1}+\cdots+s_{\tau_{2b-1}})}
G(s_1,\cdots,s_{2k}),\eqno(1.8)
$$
where $\zeta$ is the Riemann zeta-function,
and $1\le\lambda_1<\cdots<\lambda_{2a}\le 2k$,
$1\le\tau_1<\cdots<\tau_{2b-1}\le 2k$ with $a,b\ge1$. The function $G$ is
regular and bounded for $\Re s_j>-c(k)$ ($j=1,\ldots,2k$) with a constant
$c(k)>0$ which could be given explicitly.
\medskip
An appropriate shift of contours, along with I.M. Vinogradov's zero-free
region for $\zeta(s)$ (see [6, Theorem 6.1]), yields
\medskip
\noindent
{\bf Theorem.}\quad{\it  As $z$ tends to infinity, we have
$$
\eqalignno{
\e{M}_1(z)&=(1+o(1)){1\over 2\pi}\int_{-\infty}^\infty
{G_1(t)\over|\zeta(1+it)|^2}{dt\over t^2},&(1.9)\cr
\e{M}_2(z)&=(1+o(1)){3\over 4\pi}(\log z)^2\int_{-\infty}^\infty
{|\zeta(1+2it)|^2\over|\zeta(1+it)|^8}G_2(t){dt\over t^4},
&(1.10)\cr
}
$$
where $G_k(t)=G(s_1,\cdots,s_{2k})$ with $s_1,\ldots,s_k=it$,
$s_{k+1},\ldots,s_{2k}=-it$.
}
\medskip
\noindent
We note that $G_k(t)>0$. The formula $(1.9)$ is proved in [2] with an argument
different from ours. The formula $(1.10)$ is a corrected version of the relevant
claim made in [10]; there was an error in the computation of certain residues.
The advantage of our argument over that of [2] is perhaps in that ours can
give rise to $(1.10)$. Our argument should work, in principle, for any $k$.
However, the mode of shifts of contours and the arrangement of residues become
formidably complicated for $k\ge3$. Thus the general case will probably
require a new approach, though the case $k=3$ appears to be still
manageable as a direct extension of the present work. 
\bigskip
\noindent
{\bf 2.\ Proof of Theorem.}\quad We shall deal with the case $k=2$ only, for
the case $k=1$ is analogous and in fact by far simpler. Also we shall assume
that $z$ is half a large odd integer. Obviously this will make no difference.
\medskip
To begin with, let $\alpha=(\log z)^{-1}$ and $T=z^2$. Then we have
$$
\eqalignno{
\e{M}_2(z)={1\over(2\pi i)^{4}}&\int_{10\alpha-10Ti}^{10\alpha+10Ti}
\int_{4\alpha-4Ti}^{4\alpha+4Ti}\int_{2\alpha-2Ti}^{2\alpha+2Ti}
\int_{\alpha-Ti}^{\alpha+Ti}M(s_1,s_2,s_3,s_4)\cr
&\times z^{s_1+s_2+s_3+s_4}
{ds_1ds_2ds_3ds_4\over s_1s_2s_3s_4}
+O\left({(\log z)^7\over z}\right),&(2.1)\cr
}
$$
where the implied constant is absolute.
To show this, we note first that for any positive integer $d_1$
$$
{1\over2\pi i}\int_{\alpha-Ti}^{\alpha+Ti}\left({z\over
d_1}\right)^{s_1}{ds_1\over s_1}=\delta(d_1)+O\left({z\over
Td_1^\alpha}\right),\eqno(2.2)
$$
where $\delta(d)=1$ if $d<z$ and $0$ otherwise. This is of course a crude
consequence of Perron's inversion formula. Multiply both sides by the factor
$(z/d_2)^{s_2}/s_2$ with an integer $d_2>0$ and integrate with respect to $s_2$
as indicated by $(2.1)$. We have
$$
\eqalignno{
{1\over(2\pi i)^2}&\int_{2\alpha-2Ti}^{2\alpha+2Ti}
\int_{\alpha-Ti}^{\alpha+Ti}{z^{s_1+s_2}\over d_1^{s_1}
d_2^{s_2}}{ds_1ds_2\over s_1s_2}\cr
&=\delta(d_1)\delta(d_2)+O\left({z\delta(d_1)\over
Td_2^\alpha}\right)
+O\left({z\log T\over
T(d_1d_2)^\alpha}\right)\cr &=\delta(d_1)\delta(d_2)+O\left({z\log T\over
T(d_1d_2)^\alpha}\right).
&(2.3)}
$$
Repeating the same procedure, we get
$$
\eqalignno{
{1\over(2\pi i)^4}&\int_{10\alpha-10Ti}^{10\alpha+10Ti}
\int_{4\alpha-4Ti}^{4\alpha+4Ti}\int_{2\alpha-2Ti}^{2\alpha+2Ti}
\int_{\alpha-Ti}^{\alpha+Ti}
{z^{s_1+s_2+s_3+s_4}\over d_1^{s_1}d_2^{s_2}d_3^{s_3}d_4^{s_4}}
{ds_1ds_2ds_3ds_4\over s_1s_2s_3s_4}\cr
&=\prod_{j=1}^4\delta(d_j)+O\left({z(\log T)^3\over
T(d_1d_2d_3d_4)^\alpha}\right).
&(2.4)
}
$$
Then, we divide both sides by $[d_1,d_2,d_3,d_4]$ and sum the result. We find
that the first term on the right of $(2.1)$ is equal to
$$
\eqalignno{
&\e{M}_2(z)+O\left({z(\log T)^3\over T}\sum_{d_1=1}^\infty
\sum_{d_2=1}^\infty\sum_{d_3=1}^\infty
\sum_{d_4=1}^\infty{|\mu(d_1)\mu(d_2)\mu(d_3)\mu(d_4)|\over
[d_1,d_2,d_3,d_4](d_1d_2d_3d_4)^\alpha}\right)\cr
=&\e{M}_2(z)+O\left({z(\log T)^3\over T}\prod_p\left(1-{1\over p}
+{1\over p}\left(1+{1\over p^\alpha}\right)^4\right)\right).&(2.5)
}
$$
Observing that this Euler product is $O(\zeta^4(1+\alpha))$, we end the proof
of $(2.1)$.
\medskip
Now, let $\beta=(\log z)^{-3/4}$. We shift the contour for the
$s_4$-integral to the vertical segment $[-\beta-10Ti,-\beta+10Ti]$. In view of
$(1.8)$, we encounter poles at $s_4=-s_1, -s_2, -s_3$, and
$-(s_1+s_2+s_3)$. Computing respective residues, we have
$$
\e{M}_2(z)=\left\{\e{M}_2^{(1)}+\e{M}_2^{(2)}+\e{M}_2^{(3)}
+\e{M}_2^{(4)}\right\}(z)+O\left(z^{-\beta/2}\right).\eqno(2.6)
$$ 
by virtue of Vinogradov's zero-free region for $\zeta$ together 
with the related bounds for $\zeta$, $1/\zeta$ (see [6, Theorem 6.3; Lemma
12.3]); the same combination will be implicitly invoked in what follows as well.
Here
$$
\eqalignno{
\e{M}_2^{(1)}(z)=&-{1\over(2\pi i)^{3}}\int_{4\alpha-4Ti}^{4\alpha+4Ti}
\int_{2\alpha-2Ti}^{2\alpha+2Ti}
\int_{\alpha-Ti}^{\alpha+Ti}
{\zeta(1+s_1+s_2)\zeta(1+s_1+s_3)
\over\zeta(1+s_1)\zeta(1+s_2)^2\zeta(1+s_3)^2}\cr
&\times{\zeta(1+s_2+s_3)^2\zeta(1+s_2-s_1)
\zeta(1+s_3-s_1)\over\zeta(1-s_1)
\zeta(1+s_1+s_2+s_3)\zeta(1+s_2+s_3-s_1)}\cr
&\times G(s_1,s_2,s_3,-s_1)z^{s_2+s_3}{ds_1ds_2ds_3\over s_1^2s_2s_3},&(2.7)\cr
\e{M}_2^{(2)}(z)=&-{1\over(2\pi i)^{3}}
\int_{4\alpha-4Ti}^{4\alpha+4Ti}
\int_{2\alpha-2Ti}^{2\alpha+2Ti}
\int_{\alpha-Ti}^{\alpha+Ti}
{\zeta(1+s_1+s_2)\zeta(1+s_1+s_3)^2
\over\zeta(1+s_1)^2\zeta(1+s_2)\zeta(1+s_3)^2},\cr
&\times{\zeta(1+s_2+s_3)\zeta(1+s_1-s_2)
\zeta(1+s_3-s_2)\over\zeta(1-s_2)
\zeta(1+s_1+s_2+s_3)\zeta(1+s_1+s_3-s_2)}\cr
&\times G(s_1,s_2,s_3,-s_2)z^{s_1+s_3}{ds_1ds_2ds_3\over s_1s_2^2s_3},&(2.8)\cr
\e{M}_2^{(3)}(z)=&-{1\over(2\pi i)^{3}}
\int_{4\alpha-4Ti}^{4\alpha+4Ti}
\int_{2\alpha-2Ti}^{2\alpha+2Ti}
\int_{\alpha-Ti}^{\alpha+Ti}
{\zeta(1+s_1+s_2)^2\zeta(1+s_1+s_3)
\over\zeta(1+s_1)^2\zeta(1+s_2)^2\zeta(1+s_3)}\cr
&\times{\zeta(1+s_2+s_3)\zeta(1+s_1-s_3)
\zeta(1+s_2-s_3)\over\zeta(1-s_3)
\zeta(1+s_1+s_2+s_3)\zeta(1+s_1+s_2-s_3)}\cr
&\times G(s_1,s_2,s_3,-s_3)z^{s_1+s_2}{ds_1ds_2ds_3\over s_1s_2s_3^2},&(2.9)\cr
\e{M}_2^{(4)}(z)=&-{1\over(2\pi i)^{3}}
\int_{4\alpha-4Ti}^{4\alpha+4Ti}
\int_{2\alpha-2Ti}^{2\alpha+2Ti}
\int_{\alpha-Ti}^{\alpha+Ti}
{\zeta(1+s_1+s_2)\zeta(1-s_1-s_2)
\over\zeta(1+s_1)\zeta(1-s_1)
\zeta(1+s_2)}\cr
&\times{\zeta(1+s_1+s_3)\zeta(1-s_1-s_3)
\zeta(1+s_2+s_3)\zeta(1-s_2-s_3)\over\zeta(1-s_2)\zeta(1+s_3)\zeta(1-s_3)
\zeta(1+s_1+s_2+s_3)
\zeta(1-s_1-s_2-s_3)}\cr
&\times G(s_1,s_2,s_3,-s_1-s_2-s_3)
{ds_1ds_2ds_3\over s_1s_2s_3(s_1+s_2+s_3)}.&(2.10)
}
$$
\medskip
Let us first show that
$$
\e{M}_2^{(4)}(z)\ll(\log z)^{3/2}.\eqno(2.11)
$$
We note that the bound $\e{M}_2^{(4)}(z)\ll 1$ appears highly probable; in
fact, this holds under the Riemann Hypothesis. To prove $(2.11)$ we
observe first that
$$
\e{M}_2^{(4)}(z)=
-{1\over(2\pi i)^{3}}
\int_{\beta-4Ti}^{\beta+4Ti}
\int_{\beta-2Ti}^{\beta+2Ti}
\int_{\beta-Ti}^{\beta+Ti}
\{\cdots\}
{ds_1ds_2ds_3\over s_1s_2s_3(s_1+s_2+s_3)}+o(1).
\eqno(2.12)
$$
Then we shift the contour of the inner-most integral to
$$
C=\left\{s_1:\,{c\over(\log(2+|t|+|s_2|+|s_3|))^{3/4}}+it,\quad -T\le t\le T
\right\},\eqno(2.13)
$$
where $c>0$ needs to be sufficiently small. We have
$$
\e{M}_2^{(4)}(z)=
-{1\over(2\pi i)^{3}}
\int_{\beta-4Ti}^{\beta+4Ti}
\int_{\beta-2Ti}^{\beta+2Ti}
\int_{C}
\{\cdots\}
{ds_1ds_2ds_3\over s_1s_2s_3(s_1+s_2+s_3)}+o(1).
\eqno(2.14)
$$
This implies that
$$
\eqalignno{
\e{M}_2^{(4)}(z)\ll &(\log z)^{3/2}\int_{-4T}^{4T}\int_{-2T}^{2T}
\int_{-T}^T\log^{10}(2+|t_1|+|t_2|+|t_3|)\cr
&\times{dt_1dt_2dt_3\over
(1+|t_1|)(1+|t_2|)(1+|t_3|)(1+
|t_1+t_2+t_3|)}+o(1).&(2.15)
}
$$
On the right side, the factor $(\log z)^{3/2}$ comes from
the factor $\zeta(1+s_2+s_3)\zeta(1-s_2-s_3)$ in $(2.14)$. In
the integrand, the denominator comes from that in $(2.14)$, and the logarithmic
factor from those zeta-factors there, save for
$\zeta(1+s_2+s_3)\zeta(1-s_2-s_3)$. One may see readily that
$$
\int_{-T}^T{\log^{10}(2+|t_1|+|t_2|+|t_3|)\over
(1+|t_1|)(1+|t_1+t_2+t_3|)}dt_1\ll{\log^{11}(2+|t_2|+|t_3|)
\over1+|t_2+t_3|}.\eqno(2.16)
$$
Thus we have 
$$
\e{M}_2^{(4)}(z)\ll(\log z)^{3/2}\int_{-4T}^{4T}{\log^{12}(2+|t_3|)\over
(1+|t_3|)^2}dt_3,\eqno(2.17)
$$
which proves our claim $(2.11)$.
\medskip
Let us treat $\e{M}_2^{(1)}$. In $(2.7)$ we shift the $s_3$-contour
to the segment $[-\beta-4Ti,-\beta+4Ti]$. We encounter poles at $s_3=s_1, -s_1,
-s_2$. Computing the respective residues we get
$$
\e{M}_2^{(1)}(z)=\left\{\e{M}_2^{(1,1)}+\e{M}_2^{(1,2)}+\e{M}_2^{(1,3)}\right\}
(z)+O(z^{-\beta/2}).\eqno(2.18)
$$
We have
$$
\e{M}_2^{(1,3)}(z)\ll \log z.\eqno(2.19)
$$
In fact, $\e{M}_2^{(1,3)}(z)$ is a linear polynomial in $\log z$, whose
coefficients are bounded. More precisely, the leading coefficient
is equal to
$$
\eqalignno{
{1\over(2\pi i)^{2}}
\int_{2\alpha-2Ti}^{2\alpha+2Ti}
\int_{\alpha-Ti}^{\alpha+Ti}&
{\zeta(1+s_1+s_2)\zeta(1-s_1-s_2)\zeta(1+s_1-s_2)\zeta(1-s_1+s_2)
\over(\zeta(1+s_1)\zeta(1-s_1)\zeta(1+s_2)\zeta(1-s_2))^2}\cr
&\times G(s_1,s_2,-s_2,-s_1){ds_1ds_2\over
(s_1s_2)^2}.&(2.20)\cr 
}
$$
We shift the contour of the outer integral to 
$$
\left\{s_2:\,{c\over(\log(2+|t|+|s_1|))^{3/4}}+it,\,-2T\le t\le
2T\right\},\eqno(2.21)
$$
with a small $c>0$. We do not encounter any pole. The new double integral is
bounded by a constant multiple of
$$
\int_{-2T}^{2T}\int_{-T}^T{\log^{12}(2+|t_1|+|t_2|)\over
((1+|t_1|)(1+|t_2|))^2}dt_1dt_2\ll 1,\eqno(2.22)
$$
as claimed. The constant term of the linear polynomial has more complicated
expression than $(2.20)$, involving derivatives of the zeta-function. However,
its treatment is analogous.
\par
On the other hand we have
$$
\eqalignno{
\e{M}_2^{(1,1)}(z)=&-{1\over (2\pi
i)^2}\int_{2\alpha-2Ti}^{2\alpha+2Ti}\int_{\alpha-Ti}^{\alpha+Ti}
{\zeta(1+s_1+s_2)^3\zeta(1+2s_1)\over\zeta(1+s_1)^3\zeta(1+s_2)^3}\cr
&\times{\zeta(1+s_2-s_1)
\over\zeta(1-s_1)
\zeta(1+2s_1+s_2)}G(s_1,s_2,s_1,-s_1)z^{s_1+s_2}
{ds_1ds_2\over s_1^3s_2},&(2.23)
\cr
\e{M}_2^{(1,2)}(z)=&{1\over (2\pi i)^2}
\int_{2\alpha-2Ti}^{2\alpha+2Ti}\int_{\alpha-Ti}^{\alpha+Ti}
{\zeta(1-s_1+s_2)^3\zeta(1-2s_1)\over\zeta(1-s_1)^3\zeta(1+s_2)^3}\cr
&\times{\zeta(1+s_2+s_1)
\over\zeta(1+s_1)
\zeta(1-2s_1+s_2)}G(s_1,s_2,-s_1,-s_1)z^{s_2-s_1}
{ds_1ds_2\over s_1^3s_2}.&(2.24)
}
$$
In the latter we shift the $s_1$-contour to the segment $[-\alpha-iT,
-\alpha+iT]$. We do not encounter any pole. In the new $s_1$-integral we
perform the change of variable $s_1\mapsto -s_1$. On noting
$G(-s_1,s_2,s_1,s_1) =G(s_1,s_2,s_1,-s_1)$, we have
$$
\left\{\e{M}_2^{(1,1)}+\e{M}_2^{(1,2)}\right\}(z)=2\e{M}_2^{(1,1)}(z)+
o(1).\eqno(2.25)
$$
\par
We then shift the $s_2$-contour in $(2.23)$ to the segment $[-\beta-2iT,
-\beta+2iT]$. We encounter poles at $s_2=s_1, -s_1$, with the resulting double
integral being $O(z^{-\beta/2})$. The first pole contributes
$$
-{1\over 2\pi i}\int_{\alpha-Ti}^{\alpha+Ti}{\zeta(1+2s_1)^4G(s_1,s_1,s_1,
-s_1)z^{2s_1}\over\zeta(1+s_1)^6\zeta(1-s_1)\zeta(1+3s_1)}{ds_1\over s_1^4},
\eqno(2.26)
$$
which is obviously $O(z^{-\beta/2})$. Thus, computing the residue at $s_2=-s_1$,
we have
$$
\eqalignno{
&\left\{\e{M}_2^{(1,1)}+\e{M}_2^{(1,2)}\right\}(z)\cr
&={(\log z)^2\over2\pi i}
\int_{\alpha-iT}^{\alpha+iT}{\zeta(1+2s_1)\zeta(1-2s_1)\over
(\zeta(1+s_1)\zeta(1-s_1))^4}G(s_1,-s_1,s_1,-s_1){ds_1\over s_1^4}
+O(\log z).&(2.27)
}
$$
This error term is actually equal to a negligible term plus a linear polynomial
of $\log z$, the coefficients of which are easily seen to be bounded.
\par
From $(2.18)$, $(2.19)$ and $(2.27)$ we obtain
$$
\e{M}_2^{(1)}(z)=(1+o(1)){1\over2\pi}(\log z)^2
\int_{-\infty}^\infty{|\zeta(1+2it)|^2\over|\zeta(1+it)|^8}
G(it,it,-it,-it){dt\over t^4}.
\eqno(2.28)
$$
which ends our computation of $\e{M}_2^{(1)}(z)$.
\medskip
Next, we shall consider $\e{M}_2^{(2)}$; we may be brief. In $(2.8)$ we shift
the
$s_3$-contour to $[-\beta-4iT,-\beta+4iT]$. We encounter poles at
$s_3=s_2,-s_2,-s_1$. Computing the respective residues, we have
$$
\e{M}_2^{(2)}(z)=\left\{\e{M}_2^{(2,1)}+\e{M}_2^{(2,2)}+\e{M}_2^{(2,3)}\right\}
(z)+O(z^{-\beta/2}).\eqno(2.29)
$$
We have $\e{M}_2^{(2,3)}(z)\ll\log z$ similarly to $(2.19)$. We have
$$
\eqalignno{
\e{M}_2^{(2,1)}(z)=&-{1\over(2\pi i)^2}
\int_{2\alpha-2Ti}^{2\alpha+2Ti}
\int_{\alpha-Ti}^{\alpha+Ti}
{\zeta(1+s_1+s_2)^3\zeta(1+2s_2)
\over\zeta(1+s_1)^3\zeta(1+s_2)^3},\cr
&\times{\zeta(1+s_1-s_2)\over\zeta(1-s_2)
\zeta(1+s_1+2s_2)}
G(s_1,s_2,s_2,-s_2)z^{s_1+s_2}{ds_1ds_2\over s_1s_2^3}, &(2.30)\cr
\e{M}_2^{(2,2)}(z)=&{1\over(2\pi i)^2}
\int_{2\alpha-2Ti}^{2\alpha+2Ti}
\int_{\alpha-Ti}^{\alpha+Ti}
{\zeta(1+s_1-s_2)^3\zeta(1-2s_2)
\over\zeta(1+s_1)^3\zeta(1-s_2)^3},\cr
&\times{\zeta(1+s_1+s_2)\over\zeta(1+s_2)
\zeta(1+s_1-2s_2)}
G(s_1,s_2,-s_2,-s_2)z^{s_1-s_2}{ds_1ds_2\over s_1s_2^3}.&(2.31)
}
$$
In the latter we shift the $s_2$-contour to the segment $[\beta-2iT,
\beta+2iT]$, and we get $\e{M}_2^{(2,2)}(z)\ll z^{-\beta/2}$. On the other hand,
in the former we shift the $s_2$-contour to $[-\beta-2iT,-\beta+2iT]$. We have
$$
\e{M}_2^{(2,1)}(z)
={(\log z)^2\over4\pi i}
\int_{\alpha-iT}^{\alpha+iT}{\zeta(1+2s_1)\zeta(1-2s_1)\over
(\zeta(1+s_1)\zeta(1-s_1))^4}G(s_1,-s_1,-s_1,s_1){ds_1\over s_1^4}
+O(\log z).\eqno(2.32)
$$
Hence we obtain
$$
\e{M}_2^{(2)}(z)=(1+o(1)){1\over4\pi}(\log z)^2
\int_{-\infty}^\infty{|\zeta(1+2it)|^2\over|\zeta(1+it)|^8}
G(it,it,-it,-it){dt\over t^4}.
\eqno(2.33)
$$
\medskip
It now remains for us to consider $\e{M}_2^{(3)}$. This time we shift first
the $s_3$-contour  in $(2.9)$ to the segment $[2\beta-4iT,2\beta+4iT]$. We do
not encounter any pole, and thus $\e{M}_2^{(3)}(z)$ is equal to the new integral
plus a negligible error. In the new integral we shift the $s_2$-contour to
the segment $[-\beta-2iT,-\beta+2iT]$. We encounter only one pole at $s_2=-s_1$.
Computing the residue we get
$$
\e{M}_2^{(3)}(z)\ll \log z.\eqno(2.34)
$$
\medskip
Finally, collecting $(2.6)$, $(2.11)$, $(2.28)$, $(2.33)$ and $(2.34)$,
we end our proof of $(1.10)$.
\bigskip
\noindent
{\csc Remark.} There is an old conjecture by P. Erd\"os about the size of
the arithmetic function
$$
\sup_z\Bigg|\sum_{d|n,d\le z}\mu(d)\Bigg|.\eqno(2.35)
$$
See [2] for details. Our problem is certainly related to the dual of
Erd\"os'; that is, the supremum is taken in $n$ instead of $z$. As to the
possible relation of our problem with the Selberg sieve, see [1], [8], [3] and
[7], in chronological order. In addition to these, see 
[8, \S1.3] for an extension of $(1.3)$. It should be noted that [8], [3] and [7]
were developed in conjunction with the zero-density theory for the Riemann
zeta- and Dirichlet $L$-functions. However, for this particular purpose those
works turned out later to be redundant due to the observation [8, (1.3.12)]. 
\vskip 2cm
\centerline{\bf References}
\bigskip
\item{[1]} M.B. Barban and P.P. Vehov. An extremal problem. Trans.\ Moscow
Math.\ Soc., {\bf 18} (1968), 91-99.
\item{[2]} F. Dress, H. Iwaniec and G. Tennenbaum. Sur une somme li\'ee
\`a la fonction M\"obius. J.\ reine angew.\ Math., {\bf 340} (1983), 53-58.
\item{[3]} S. Graham. An asymptotic estimate related to Selberg's sieve. J.\
Number Theory, {\bf 10} (1978), 83-94.
\item{[4]} D. Goldston and C.M. Y{\i}ld{\i}r{\i}m. Higher correlations of
divisor sums related to primes.\ III. Preprint, September 2002.
\item{[5]} ---. Small gaps between primes. Notes, March 2003.
\item{[6]} A. Ivi\'c. {\it The Riemann Zeta-Function.\ Theory and
Applications.\/} Dover Publ., Inc., Mineola, New York 2003.
\item{[7]} M. Jutila. On a problem of Barban and Vehov. Mathematika, {\bf 26}
(1979), 62-71. 
\item{[8]} Y. Motohashi. A problem in the theory of sieves. Kokyuroku RIMS
Kyoto Univ., {\bf 222} (1974), 9-50.
\item{[9]} ---. {\it Sieve Methods and Prime Number Theory.\/} Tata IFR, 
Lect.\ Math.\ Phys., {\bf 72}, Springer Verlag, Berlin etc.\ 1983.
\item{[10]} ---. M\"obius function over divisors. Notes, March 1985.

\vskip 1cm
\noindent
Yoichi Motohashi
\smallskip\noindent
Honkomagome 5-67-1 \#901, Tokyo 113-0021, Japan
\par\noindent
Email: am8y-mths@asahi-net.or.jp
\par\noindent
HP: www.ne.jp/asahi/zeta/motohashi/

\bye